\documentclass[12pt]{amsart}
\usepackage{amsfonts}
\usepackage{amsmath,amssymb,amsthm}

\newcommand{\R}{\mathbb R}

\newcommand{\E}{\mathbb E}

\renewcommand{\span}{\mathrm{span}}

\newtheorem{thm}{Theorem}[section]

\newtheorem{lem}[thm]{Lemma}
\newtheorem{prop}[thm]{Proposition}
\theoremstyle{definition}
\newtheorem{defn}[thm]{Definition}
\theoremstyle{remark}

\newcommand{\ds}{\displaystyle}

\begin{document}

\title[MINIMAL SURFACES IN THE FOUR-DIMENSIONAL EUCLIDEAN SPACE]
{MINIMAL SURFACES IN THE FOUR-DIMENSIONAL\\ EUCLIDEAN SPACE}

\author{Georgi Ganchev and Velichka Milousheva}
\address{Bulgarian Academy of Sciences, Institute of Mathematics and Informatics,
Acad. G. Bonchev Str. bl. 8, 1113 Sofia, Bulgaria}
\email{ganchev@math.bas.bg}
\address{Bulgarian Academy of Sciences, Institute of Mathematics and Informatics,
Acad. G. Bonchev Str. bl. 8, 1113, Sofia, Bulgaria}
\email{vmil@math.bas.bg}



\begin{abstract}
We prove that the Gauss curvature and the curvature of the normal connection of
any minimal surface in the four dimensional Euclidean space satisfy an inequality,
which generates two classes of minimal surfaces: minimal surfaces of general type
and minimal super-conformal surfaces. We prove a Bonnet-type theorem for strongly
regular minimal surfaces of general type in terms of their invariants.
We introduce canonical parameters on strongly regular minimal surfaces of general
type and prove that any such a surface is determined up to a motion by two invariant
functions satisfying a system of two natural partial differential equations. On any
minimal surface of the basic class of non strongly regular minimal surfaces we
define canonical parameters and prove that any such a surface is determined up to
a motion by two invariant functions of one variable satisfying a system of two
natural ordinary differential equations. We find a geometric description of this
class of non strongly regular minimal surfaces.
\end{abstract}

\maketitle

\section{Introduction} \label{S:Intr}

In \cite{GMih} it is proved that any strongly regular Weingarten surface in Euclidean
space ${\R}^3$ admits geometrically determined parameters, which are an analogue of
the natural parameter in the theory of the smooth curves in ${\R}^3$. As a corollary
of this fact it follows that any such a surface is determined uniquely up to a motion
by its normal curvature function satisfying the corresponding natural partial
differential equation.

In this paper we study minimal surfaces in the four dimensional Euclidean space $\R^4$.
As usual, the notion of a minimal surface means a regular surface with zero mean
curvature vector. At any point of a regular surface $M^2$ in ${\R}^4$ in \cite{GMil}
we introduced two invariants $k$ and $\varkappa$. In terms of these invariants minimal
surfaces are characterized by the equality
$$\varkappa^2-k=0.$$
The function $\varkappa$ is the curvature of the normal connection of $M^2$. Denoting
by $K$ the Gauss curvature of $M^2$, in Section 3 we prove that on any minimal
surface the following inequality holds
$$K^2-\varkappa^2\geq 0.$$
This inequality generates two classes of minimal surfaces:
\begin{itemize}
\item
the class of minimal super-conformal surfaces characterized by $K^2 - \varkappa^2 =0$;
\item
the class of minimal surfaces of general type characterized by $K^2-\varkappa^2>0$.
\end{itemize}

In this paper we study the class of minimal surfaces of general type.

In Section 4 we introduce semi-canonical parameters and canonical lines on any minimal
surface of general type, which play the same role as the principal parameters and the
lines of curvature in the theory of surfaces in ${\R}^3$.

In Section 5 we denote by $\gamma_1$ and $\gamma_2$ the geodesic curvatures of the
canonical lines on a minimal  surface $M^2$ of general type, and define strongly regular
minimal surfaces of general type by the inequality
$$\gamma_1\,\gamma_2\neq 0.$$
Further we prove the fundamental theorem of Bonnet-type (Theorem \ref{T:Fundamental Theorem})
for strongly regular minimal surfaces of general type in terms of four invariant functions.

Next we define canonical parameters and prove (Theorem \ref{T:Canonocal parameters})
that any strongly regular minimal surface of general type admits such parameters. The
fundamental theorem (Theorem \ref{T:Fundamental Theorem - Canonical}) in canonical
parameters states as follows:

{\it Any two solutions $\mu(u,v)$ and $\nu(u,v)$ to the system of natural partial
differential equations
$$\frac{1}{4}\sqrt{|\mu^2-\nu^2|} \, \Delta \ln |\mu^2-\nu^2|+ \mu^2+\nu^2=0,$$
$$\frac{1}{2}\sqrt{|\mu^2-\nu^2|} \, \Delta \ln \left|\frac{\mu+\nu}{\mu-\nu}\right|+ 2\mu\nu=0,$$
determine a unique (up to a motion) strongly regular minimal surface of general type
with invariants $\mu$ and $\nu$. Furthermore $(u,v)$ are canonical parameters.}

In Section 7 we consider non strongly regular minimal surfaces of general type. The
basic class of these surfaces is determined by the condition $\gamma_1 = 0$ , i.e.
one of the families of canonical lines consists of geodesics.

We introduce canonical parameters and prove (Theorem 7.3) that any minimal surface
of general type satisfying the condition $\gamma_1=0$ admits canonical parameters.
The fundamental theorem (Theorem 7.4) for these surfaces in canonical parameters
states as follows:
\vskip 2mm
{\it Any two solutions $\mu(u)$ and $\nu(u)$ to the system of natural ordinary
differential equations
$$\begin{array}{l}
\ds{\frac{1}{4}\left(\ln |\mu^2 - \nu^2|\right)_{uu}
- \frac{1}{16}\left(\ln |\mu^2 - \nu^2|\right)_u^2  + \nu^2 + \mu^2 = 0},\\
[3mm]
\ds{\frac{1}{2}\left(\ln \left|\frac{\mu + \nu}{\mu - \nu}\right|\right)_{uu}
- \frac{1}{8}\left(\ln |\mu^2 - \nu^2|\right)_u \left(\ln \left|\frac{\mu
+ \nu}{\mu - \nu} \right|\right)_u + 2\nu \mu = 0}
\end{array}$$
determine a unique (up to a motion) non strongly regular surface $(\gamma_1=0)$
with invariants $\mu$ and $\nu$.}
\vskip 2mm
We give a geometric description of non strongly regular minimal surfaces of general type
(with $\gamma_1=0$) in Proposition 7.5.

\section{Preliminaries} \label{S:Pre}

We denote by $g$ the standard metric in the four-dimensional Euclidean space
$\R^4$ and by $\nabla'$ its flat Levi-Civita
connection. All considerations in the present paper are local and
all functions, curves, surfaces, tensor fields etc. are assumed to
be of the class $\mathcal C^{\infty}$.

Let $M^2: z = z(u,v), \, \, (u,v) \in {\mathcal D}$ (${\mathcal D} \subset \R^2$)
be a regular 2-dimensional surface in $\R^4$. The tangent space to $M^2$ at
an arbitrary point $p=z(u,v)$ of $M^2$ is ${\rm span} \{z_u, z_v\}$.

For an arbitrary orthonormal normal frame field $\{e_1, e_2\}$ of
$M^2$ we have the standard derivative formulas:
$$\begin{array}{l}
\vspace{2mm} \nabla'_{z_u}z_u=z_{uu} = \Gamma_{11}^1 \, z_u +
\Gamma_{11}^2 \, z_v
+ c_{11}^1\, e_1 + c_{11}^2\, e_2,\\
\vspace{2mm} \nabla'_{z_u}z_v=z_{uv} = \Gamma_{12}^1 \, z_u +
\Gamma_{12}^2 \, z_v
+ c_{12}^1\, e_1 + c_{12}^2\, e_2,\\
\vspace{2mm} \nabla'_{z_v}z_v=z_{vv} = \Gamma_{22}^1 \, z_u +
\Gamma_{22}^2 \, z_v
+ c_{22}^1\, e_1 + c_{22}^2\, e_2,\\
\end{array}$$
where $\Gamma_{ij}^k$ are the Christoffel's symbols and $c_{ij}^k$,
$i, j, k = 1,2$ are functions on $M^2$.

We use the standard denotations \;$E(u,v)=g(z_u,z_u), \; F(u,v)=g(z_u,z_v), \;
G(u,v)=g(z_v,z_v)$ for the coefficients of the first fundamental form and set
$W=\sqrt{EG-F^2}$. Denoting by $\sigma$ the second fundamental tensor of $M^2$,
we have
$$\begin{array}{l}
\sigma(z_u,z_u)=c_{11}^1\, e_1 + c_{11}^2\, e_2,\\
[2mm]
\sigma(z_u,z_v)=c_{12}^1\, e_1 + c_{12}^2\, e_2,\\
[2mm] \sigma(z_v,z_v)=c_{22}^1\, e_1 + c_{22}^2\, e_2.\end{array}$$
\vskip 2mm

In \cite{GMil} we found invariants of a surface $M^2$ in  $\R^4$,
defining a geometrically determined linear map in the tangent space of the surface.
We introduce the functions
$$\Delta_1 = \left|%
\begin{array}{cc}
\vspace{2mm}
  c_{11}^1 & c_{12}^1 \\
  c_{11}^2 & c_{12}^2 \\
\end{array}%
\right|, \quad
\Delta_2 = \left|%
\begin{array}{cc}
\vspace{2mm}
  c_{11}^1 & c_{22}^1 \\
  c_{11}^2 & c_{22}^2 \\
\end{array}%
\right|, \quad
\Delta_3 = \left|%
\begin{array}{cc}
\vspace{2mm}
  c_{12}^1 & c_{22}^1 \\
  c_{12}^2 & c_{22}^2 \\
\end{array}%
\right|;$$
$$L(u,v) = \displaystyle{\frac{2 \Delta_1}{W}, \quad M(u,v) =
\frac{\Delta_2}{W}, \quad N(u,v) = \frac{2 \Delta_3}{W}}.$$
Further we denote
$$\begin{array}{l}
\vspace{2mm}
\displaystyle{\gamma_1^1=\frac{FM-GL}{EG-F^2}, \qquad \gamma_1^2
=\frac{FL-EM}{EG-F^2}},\\
\vspace{2mm}
\displaystyle{\gamma_2^1=\frac{FN-GM}{EG-F^2}, \qquad
\gamma_2^2=\frac{FM-EN}{EG-F^2}},
\end{array}$$
and consider the linear map
$$\gamma: T_pM^2 \rightarrow T_pM^2,$$
determined by the conditions
$$\begin{array}{l}
\vspace{2mm}
\gamma(z_u)=\gamma_1^1z_u+\gamma_1^2z_v,\\
\vspace{2mm} \gamma(z_v)=\gamma_2^1z_u+\gamma_2^2z_v.
\end{array}$$
The linear map $\gamma$ of Weingarten type at the point $p \in M^2$ is invariant
up to a sign and the functions
$$k = \frac{LN - M^2}{EG - F^2}, \qquad
\varkappa =\frac{EN+GL-2FM}{2(EG-F^2)}$$ are invariants of the surface $M^2$.
\vskip 2mm
The invariants $k$ and $\varkappa$ divide the points of $M^2$ into four types:
flat, elliptic, parabolic and hyperbolic. The surfaces consisting of flat points
satisfy the conditions
$$k(u,v)=0, \quad \varkappa(u,v)=0, \qquad (u,v) \in \mathcal D,$$
or equivalently
$$L(u,v)=0, \quad M(u,v)=0, \quad N(u,v)=0, \qquad (u,v) \in \mathcal D.$$
They are either planar surfaces (there exists a hyperplane $\R^3 \subset \R^4$
containing $M^2$) or developable ruled surfaces \cite{GMil}.
\vskip 2mm
Any surface $M^2$ in $\R^4$ satisfies the following inequality:
$$\varkappa^2 - k \geq 0.$$

The minimal surfaces in $\R^4$ are characterized by

\begin{prop}\label{P:minimal} \cite{GMil}
A surface $M^2$ in $\R^4$ is minimal if and only if
$$\varkappa^2 - k = 0.$$
\end{prop}

\section{Classes of minimal surfaces}\label{S:Minimal}

The minimal surfaces in $\R^4$ can also be characterized in terms of the ellipse
of curvature. Let us recall that the \textit{ellipse of curvature} of a surface
$M^2$ in $\R^4$ at a point $p \in M^2$ is the ellipse in the normal space of $M^2$
at the point $p$ given by $\{\sigma(x,x): \, x \in T_pM^2, \, g(x,x) = 1\}$
\cite{MW1, MW2}. Let $\{x,y\}$ be an orthonormal base of the tangent space
$T_pM^2$ at $p$. Then, we have for any $v = \cos \theta \, x + \sin \theta \, y$ that
$$\sigma(v, v) = H + \ds{\cos 2\theta \, \frac{\sigma(x,x) - \sigma(y,y)}{2}
+ \sin 2 \theta\, \sigma(x,y)},$$
where $H = \ds{\frac{1}{2} \left(\sigma(x,x) + \sigma(y,y)\right)}$ is the mean
curvature normal vector of $M^2$ at $p$. So, when $v$ goes once around the unit
tangent circle, the vector $\sigma(v,v)$ goes twice around an ellipse centered
at $H$. The vectors $\ds{\frac{\sigma(x,x) - \sigma(y,y)}{2}}$ \, and $\sigma(x,y)$
determine conjugate directions of the ellipse.

A surface $M^2$ in $\R^4$ is called \textit{super-conformal} \cite{DT} if at
any point of $M^2$ the ellipse of curvature is a circle.
In \cite{DT} it is given an explicit construction of any simply connected
super-conformal surface in $\R^4$ that is free of minimal and flat points.
\vskip 2mm
Obviously, a surface $M^2$ is minimal, i.e. the mean curvature vector field $H = 0$,
if and only if for each point $p \in M^2$ the ellipse of curvature is centered at
the point $p$.

Hence, the following proposition holds.

\begin{prop}\label{P:minimal}
Let $M^2$ be a surface in $\R^4$. Then the following conditions are equivalent:

(i) $M^2$ is minimal;

(ii) $\varkappa^2 - k = 0$;

(iii) for each point $p \in M^2$ the ellipse of curvature is centered at $p$.
\end{prop}
\vskip 3mm
Now we shall consider a minimal surface $M^2$ in $\R^4$ . Without
loss of generality we assume that
$\{\displaystyle{x=\frac{z_u}{\sqrt E}, \; y=\frac{z_v}{\sqrt G}}\}$
is an orthonormal tangent frame field, i.e. $F = 0$. Let $\{e_1, e_2\}$ be a normal
frame field of  $M^2$, such that $\{x, y, e_1, e_2\}$ is a positive oriented
orthonormal frame field in $\R^4$. Since $M^2$ is minimal, then
$\sigma(x,x) + \sigma(y,y) = 0$. So with respect to the frame field $\{x, y, e_1, e_2\}$
the derivative formulas of $M^2$ get the form:
$$\begin{array}{ll}
\vspace{2mm}
\nabla'_xx = \quad \quad \quad \gamma_1\,y + \,a\, e_1 + b\,e_2,  & \qquad
\nabla'_x e_1 = - a\,x - c\,y \quad \quad \quad \,\, + \beta_1\,e_2,\\
\vspace{2mm}
\nabla'_xy = - \gamma_1\,x \quad \quad \; + \; c\, e_1 + d\, e_2,  & \qquad
\nabla'_y e_1 = - c \,x + \; a \,y \quad \quad \quad + \beta_2\, e_2,\\
\vspace{2mm}
\nabla'_yx = \quad \quad \; - \gamma_2\,y \; + c\, e_1 + d\, e_2,  & \qquad
\nabla'_x e_2= - b \,x - d\,y - \beta_1\, e_1,\\
\vspace{2mm}
\nabla'_yy = \;\;\gamma_2\,x \quad \quad \,\,\, - \,a\,
e_1 - b\,e_2, & \qquad \nabla'_y e_2 = - d\,x + b \,y -\beta_2\,e_1,
\end{array} \leqno{(3.1)}$$
where $\gamma_1 = - y(\ln \sqrt{E}), \,\, \gamma_2 = - x(\ln \sqrt{G})$, $a$, $b$,
$c$, $d$ are functions on $M^2$.

The formulas (3.1) imply that
$$L = 2(a d - b c)\,E, \qquad  M = 0, \qquad N = 2(a d - b c)\,G.$$
In case $L = M = N = 0$ the surface $M^2$ consists only of flat points and is
contained in 3-dimensional space $\R^3$, i.e. $M^2$ is a minimal surface in $\R^3$.
We shall consider only the non-trivial case $(L,M,N)\neq ( 0,0,0)$. So, we assume
that $a d - b c \neq 0$. The derivative formulas (3.1) imply that the invariants
$k$ and $\varkappa$ of $M^2$ are expressed as follows:
$$k = 4(a d - b c)^2, \qquad \varkappa = 2 (a d - b c). \leqno{(3.2)}$$

According to the Gauss equation the Gauss curvature $K$ of $M^2$ is expressed by
$K = g(\sigma(x,x), \sigma(y,y)) - g(\sigma (x,y),\sigma(x,y))$, and hence
$$K = - (a^2 + b^2 + c^2 + d^2). \leqno{(3.3)}$$
\noindent
\textbf{Remark:} Obviously, the Gauss curvature of each minimal surface $M^2$ in
$\R^4$ is non-positive. The inequality $K \leq 0$ for a minimal surface $M^2$ in $\R^4$
also follows from the inequality \cite{W} $K + |\varkappa| \leq \Vert H \Vert ^2$,
which holds for an arbitrary surface in $\R^4$.
\vskip 2mm
\begin{lem}\label{L:inequality}
For each minimal surface $M^2$ in $\R^4$ the following inequality holds:
$$K^2 - \varkappa^2 \geq 0.$$
\end{lem}
\vskip 2mm
\noindent
\emph{Proof:} From (3.2)  and (3.3) we get that
$$K^2 - \varkappa^2 = (a^2 + b^2 - c^2 - d^2)^2 + 4(ac + bd)^2,$$
which implies the inequality $K^2 - \varkappa^2 \geq 0$.
\qed
\vskip 2mm
Obviously, $K^2 - \varkappa^2 = 0$ if and only if
$$a^2 + b^2 = c^2 + d^2, \qquad ac + bd = 0, $$
i.e.
$$\sigma(x,x)\,\, \bot \,\, \sigma(x,y), \qquad \sigma^2(x,x) = \sigma^2(x,y),$$
which is equivalent to the condition that the ellipse of curvature is a circle.

Hence, the class of minimal surfaces, characterized by the condition
$K^2 - \varkappa^2 = 0$, is the class of minimal super-conformal surfaces.
\vskip 2mm
We shall consider the class of minimal surfaces satisfying the condition
$K^2 - \varkappa^2 > 0$, and we shall call them \emph{minimal surfaces of general type}.

Now let $M^2: z=z(u,v), \; (u,v) \in {\mathcal D}$ be a minimal surface of general
type and $p=z(u_0,v_0)$ be a fixed point.

\begin{defn}
A  tangent at a point $p \in M^2$ generated by a unit vector $x$ is said to be
\textit{canonical} if it is collinear with an axis of the ellipse of curvature at $p$.
\end{defn}

We shall find a local orthonormal tangent frame field $\{x, y\}$, defined in
${\mathcal D}_0 \subset {\mathcal D}$ ($(u_0, v_0) \in {\mathcal D}_0$), such that
$x$ and $y$ are canonical tangents, i.e. $\sigma(x,x)\, \bot \, \sigma(x,y)$, or
equivalently $a\,c + b\,d = 0$.

If $\{\overline{x},\overline{y}\}$ is an orthonormal tangent frame field, such that
$$\begin{array}{l}
\vspace{2mm}
\overline{x} = \cos \varphi \, x + \sin \varphi \, y,\\
\vspace{2mm}
\overline{y} = - \sin \varphi \, x +  \cos \varphi \, y,
\end{array} \qquad \varphi = \angle (x, \overline{x}),$$
then the corresponding functions $\overline{a}, \, \overline{b}, \,
\overline{c}, \, \overline{d}$ are expressed as follows:
$$\begin{array}{ll}
\vspace{2mm}
\overline{a} = a\,\cos 2\varphi + c \sin 2\varphi, & \qquad
\overline{b} = b\,\cos 2\varphi + d \sin 2\varphi,\\
\vspace{2mm}
\overline{c} = c\,\cos 2\varphi - a \sin 2\varphi,  & \qquad
\overline{d} = d\,\cos 2\varphi - b \sin 2\varphi.
\end{array}$$
Hence,
$$\begin{array}{l}
\vspace{2mm} \overline{a}\, \overline{c } + \overline{b}\, \overline{d}
= (a\,c + b\,d)\, \cos 4 \varphi - \ds{\frac{a^2 + b^2 - c^2 - d^2}{2}\,
\sin 4 \varphi},\\
\vspace{2mm}
\ds{\frac{\overline{a}^2 + \overline{b}^2 -
\overline{c}^2 - \overline{d}^2}{2} = \frac{a^2 + b^2 - c^2 -
d^2}{2} \, \cos 4 \varphi + (a\,c + b\,d)\, \sin 4 \varphi}.
\end{array}$$

If $a\,c + b\,d \neq 0$ at the point $(u_0, v_0) \in {\mathcal D}$, then there
exists ${\mathcal D}_0 \subset {\mathcal D}$, such that $a\,c + b\,d \neq 0$
for all $(u, v) \in {\mathcal D}_0$. Hence, in ${\mathcal D}_0$ we can change
the tangent frame field $\{x,y\}$ with $\{\overline{x},\overline{y}\}$, where
$\ds{\cot 4 \varphi = \frac{a^2 + b^2 - c^2 - d^2}{2(a\,c + b\,d)}}$. Then we
get $\overline{a}\,\overline{c} + \overline{b}\,\overline{d} = 0$, i.e.
$\sigma(\overline{x},\overline{x})\, \bot \, \sigma(\overline{x},\overline{y})$.
Moreover, $\overline{a}^2 + \overline{b}^2 - \overline{c}^2 - \overline{d}^2 \neq 0$,
i.e.
$\sigma^2(\overline{x},\overline{x}) \neq \sigma^2(\overline{x},\overline{y})$.

If $a\,c + b\,d = 0$ at the point $(u_0, v_0) \in {\mathcal D}$,
then $a^2 + b^2 - c^2 - d^2 \neq 0$ at $(u_0, v_0)$ and there
exists ${\mathcal D}_0 \subset {\mathcal D}$, such that $a^2 + b^2 - c^2 - d^2 \neq 0$
for all $(u, v) \in {\mathcal D}_0$. So, in ${\mathcal D}_0$ we can change the tangent
frame field $\{x,y\}$ with $\{\overline{x},\overline{y}\}$, where
$\ds{\tan 4 \varphi = \frac{2(a\,c + b\,d)}{a^2 + b^2 - c^2 - d^2}}$. Then we get
$\overline{a}\,\overline{c} + \overline{b}\,\overline{d} = 0$,
i.e. $\overline{x}$ and $\overline{y}$ are canonical tangents.

The canonical tangents  are uniquely determined at any point of a minimal surface
of general type.

\vskip 2mm
So, we can find a local orthonormal normal frame field $\{n_1, n_2\}$, such that
$\{x, y, n_1, n_2\}$ is a positive oriented orthonormal frame field and the
derivative formulas of $M^2$ are as follows:
$$\begin{array}{ll}
\vspace{2mm}
\nabla'_xx = \quad \quad \quad \gamma_1\,y + \,\nu\, n_1,  & \qquad
\nabla'_x n_1 = - \nu\,x \quad \quad \quad \qquad \,\, + \beta_1\,n_2,\\
\vspace{2mm}
\nabla'_xy = - \gamma_1\,x \quad \quad \; \qquad \quad +  \mu\, n_2, & \qquad
\nabla'_y n_1 = \quad \qquad \nu \,y \quad \quad \quad + \beta_2\, n_2,\\
\vspace{2mm}
\nabla'_yx = \quad \quad \; - \gamma_2\,y \; \qquad \quad + \mu \, n_2,  & \qquad
\nabla'_x n_2= \qquad - \,\, \mu \,y - \beta_1\, n_1,\\
\vspace{2mm}
\nabla'_yy = \;\;\gamma_2\,x \quad \quad \,\,\, - \,\nu \, n_1, & \qquad \nabla'_y n_2
= - \mu \,x \qquad \,\, - \beta_2\,n_1,
\end{array} \leqno{(3.4)}$$
where $ \mu > 0,\,\, \nu \neq 0$, and  $\mu^2 \neq \nu ^2$.

The invariants $k$, $\varkappa$ and the Gauss curvature $K$ of
$M^2$ are expressed as follows:
$$ k = 4 \nu^2 \mu^2; \qquad \varkappa = 2\, \nu \mu; \qquad K = - (\nu^2 + \mu^2).$$
Consequently, on a minimal surface of general type $M^2$ in $\R^4$
there exists locally a geometric frame field $\{x,y, n_1, n_2\}$
of $M^2$, such that the Frenet type derivative formulas (3.4) hold.
\vskip 2mm
We shall use the following terminology:

The integral lines of the canonical tangent vector fields $x$ and $y$ are said to be
the \emph{canonical lines}; any parameters $(u,v)$ of $M^2$ generating
canonical parametric lines are said to be \textit{semi-canonical parameters}.

It is clear that semi-canonical parameters are determined up to
changes: $ u=u(\bar u), \,\, v=v(\bar v)$ or $ u=u(\bar v), \,\, v=v(\bar u)$.

The functions $\nu,\, \mu,\, \gamma_1,\, \gamma_2,\, \beta_1, \, \beta_2$ in the
Frenet type formulas (3.4) are invariants of the surface $M^2$. The invariants
$\nu$ and $\mu$ are determined by the canonical tangents $x$ and $y$ as follows:
$$\nu^2 = \sigma^2(x,x), \qquad \mu^2 = \sigma^2(x,y).$$
The functions $\nu^2$ and $\mu^2$ are expressed by the Gauss curvature $K$
and the invariant $\varkappa$ in the following way:
$$\nu^2 + \mu^2 = -K, \qquad 4 \nu^2\mu^2 = \varkappa^2.$$
Hence, there are two classes of minimal surfaces of general type:
the class of those ones for which
$\mu^2 - \nu^2 > 0$ (equivalently $\sigma^2(x,y) > \sigma^2(x,x)$),
and the class of those ones for which
$\mu^2 - \nu^2 < 0$ (equivalently $\sigma^2(x,y) < \sigma^2(x,x)$).
\vskip 2mm
Taking into account (3.4), the Gauss and  Codazzi equations imply the following
equalities for a minimal surface $M^2$ of general type:
$$\begin{array}{l}
\vspace{2mm}
2\mu \, \gamma_2 + \nu \,\beta_2 =  x(\mu),\\
\vspace{2mm}
2 \mu \, \gamma_1 - \nu\,\beta_1 =  y(\mu),\\
\vspace{2mm}
2\nu \, \gamma_2 + \mu\,\beta_2 = x(\nu),\\
\vspace{2mm}
2\nu \, \gamma_1 -  \mu\,\beta_1 =  y(\nu),\\
\vspace{2mm}
\gamma_2 \beta_2 - \gamma_1 \beta_1 - 2 \nu \mu = x(\beta_2) - y(\beta_1),\\
\vspace{2mm} \gamma_1^2 + \gamma_2^2 - (\nu^2 + \mu^2) =
x(\gamma_2) + y(\gamma_1).
\end{array} \leqno{(3.5)}$$
From (3.5) and $\nu \mu \neq 0$, it follows that $\gamma_1^2 +
\gamma_2^2 \neq 0$, \, $\beta_1^2 + \beta_2^2 \neq 0$.

In \cite{GMil} we proved that the inavariant $\varkappa$ of a non-minimal surface
is equal to the curvature of its normal connection. Next we prove that the same
result is also true for minimal surfaces.

\begin{prop}\label{P:normal curvature}
The invariant $\varkappa$ of a minimal surface $M^2$ is equal  to the curvature of
the normal connection of the surface.
\end{prop}

\noindent
{\it Proof:} Let $M^2$ be a minimal surface of general type. Then the derivative
formulas (3.4) imply that the curvature tensor $R^{\bot}$ of the normal connection
$D$ of $M^2$ is expressed as follows:
$$\begin{array}{l}
\vspace{2mm}
R^{\bot}_{n_1}(x,y) = D_xD_yn_1 - D_yD_xn_1 - D_{[x,y]}n_1 =
\left(x(\beta_2) - y(\beta_1) + \gamma_1\,\beta_1 - \gamma_2 \,\beta_2 \right)\,n_2,\\
\vspace{2mm}
R^{\bot}_{n_2}(x,y) = D_xD_yn_2 - D_yD_xn_2 - D_{[x,y]}ln_2
= - \left(x(\beta_2) - y(\beta_1) + \gamma_1\,\beta_1 - \gamma_2 \,\beta_2 \right)\,n_1.
\end{array}\leqno{(3.6)}$$
Taking in mind (3.5) and (3.6) we get:
$$\begin{array}{l}
\vspace{2mm}
R^{\bot}_{n_1}(x,y) = - \varkappa \,n_2,\\
\vspace{2mm} R^{\bot}_{n_2}(x,y) = \varkappa \,n_1,
\end{array}$$
i.e.
$$\varkappa =  g(R^{\bot}_{n_2}(x,y), n_1) = g(R^{\bot}(x,y) n_2, n_1). \leqno{(3.7)}$$

If $M^2$ is a minimal super-conformal surface, the formula (3.7)
can be derived in a similar way.

The function $g(R^{\bot}(x,y)n_2, n_1)$ is the curvature of the normal connection of
$M^2$. Hence, the invariant $\varkappa$ of a minimal surface $M^2$ is equal to
the curvature of the normal connection of the surface.
\qed

\section{Geometric characterization of the canonical tangents}
\label{S:Geometric characterization}

Let $M^2: z = z(u,v), \, \, (u,v) \in {\mathcal D}$ (${\mathcal D} \subset \R^2$)
be a minimal surface of general type, parameterized by semi-canonical parameters.
We shall give a geometric characterization of the canonical tangents in terms of the
geodesic curves, determined by these tangents.

\vskip 2mm
\begin{thm}\label{T:geometric characterization}
Let $M^2$ be a minimal surface of general type, parameterized by semi-canonical
parameters. Then a tangent $g$ at a point $p \in M^2$ is canonical if and only if
the principal normal of the geodesic curve, passing through the point $p$ and
tangent to $g$ at $p$, is collinear with  the geometric normal vector field
$n_1$ of $M^2$.
\end{thm}
\vskip 2mm
\noindent
\emph{Proof:} Suppose that $c: u = u(s), \,\, v = v(s), \,\, s \in J$ ($J \subset \R$)
is a curve on $M^2$, parameterized by the arc-length. Then the tangent vector $t$
of $c$ is expressed as follows:
$$t = z' = u'\, z_u + v' \,z_v,$$
and hence
$$t' = z'' = u''\, z_u + v'' \, z_v + u'^2\, z_{uu} + 2 u' v'\, z_{uv} +  v'^2\, z_{vv}.
\leqno{(4.1)}$$
From equalities (3.4) we get
$$\begin{array}{l}
\vspace{2mm}
z_{uu} = \Gamma_{11}^1 \, z_u + \Gamma_{11}^2 \, z_v + E \nu\, n_1,\\
\vspace{2mm}
z_{uv} = \Gamma_{12}^1 \, z_u + \Gamma_{12}^2 \, z_v  \qquad
\quad  + \sqrt{E G} \mu\, n_2,\\
\vspace{2mm}
z_{vv} = \Gamma_{22}^1 \, z_u + \Gamma_{22}^2 \, z_v - G \nu\, n_1,\\
\end{array} \leqno{(4.2)}$$
where $\Gamma_{ij}^k$ are the Christoffel's symbols. So, equalities
(4.1) and (4.2) imply that
$$\begin{array}{ll}
t' = & (u'' +  \Gamma_{11}^1 u'^2 + 2\Gamma_{12}^1 u' v' + \Gamma_{22}^1 v'^2)\, z_u \\
[2mm]
& + (v'' +  \Gamma_{11}^2 u'^2 + 2\Gamma_{12}^2 u' v' + \Gamma_{22}^2 v'^2)\, z_v \\
[2mm]
& + \nu (E u'^2 - G v'^2)\, n_1 + 2  \mu  \sqrt{EG}
\,u' v'\,n_2.
\end{array}$$
If in addition $c$ is a geodesic curve on $M^2$, i.e.
$$\begin{array}{l}
\vspace{2mm}
u'' +  \Gamma_{11}^1 u'^2 + 2\Gamma_{12}^1 u' v' + \Gamma_{22}^1 v'^2 = 0,\\
\vspace{2mm}
v'' +  \Gamma_{11}^2 u'^2 + 2\Gamma_{12}^2 u' v' +
\Gamma_{22}^2 v'^2 = 0,
\end{array}$$
then
$$t' = \nu (E u'^2 - G v'^2)\, n_1 + 2 \mu \sqrt{EG}\, u' v' \,n_2. \leqno{(4.3)}$$
\vskip 2mm
I. Now let us denote by $c_1$ the geodesic curve tangent to the canonical direction
$z_u$. Then its tangent vector is $t_{c_1} =  x = \ds{\frac{z_u}{\sqrt{E}}}$, i.e.
$u' = \ds{\frac{1}{\sqrt{E}}}, \,\, v' = 0$. Hence, from (4.3) we obtain
$$t_{c_1}' = \nu \, n_1.$$
Analogously, if $c_2$ is the geodesic curve tangent to the canonical
direction $z_v$, then its tangent vector $t_{c_2}$ satisfies
$$t_{c_2}' = - \nu \, n_1.$$
Consequently, the principal normals of the geodesic curves $c_1$ and
$c_2$ are collinear to the geometric normal vector field $n_1$.
\vskip 2mm
II. Now let $c$ be a geodesic curve on $M^2$ with tangent vector
$t = u'\, z_u + v' \,z_v$, such that $t'$ is collinear to the geometric normal
vector field $n_1$. Since $\mu \sqrt{EG} \neq 0$, then from (4.3) it follows that
$u' v' = 0$, i.e. $u = const$ or $v = const$.
Hence, the tangent to the geodesic curve $c$ at the point $p$ is canonical.
\qed
\vskip 2mm
So, for each minimal surface of general type there exists a geometric normal
vector field $n_1$, which is geometrically determined by the condition that
it is collinear with the principal normals of the geodesic curves tangent to
the canonical tangents.

\section{Strongly regular minimal surfaces of general type}

Let $M^2$ be a minimal surface of general type in $\R^4$, parameterized by
semi-canonical parameters. Then the equalities (3.5) hold. Using that
$x = \displaystyle{\frac{z_u}{\sqrt{E}}, \, y = \frac{z_v}{\sqrt{G}}}$, we can
rewrite  (3.5) in the following way:
$$\begin{array}{l}
\vspace{2mm}
2 \mu \, \gamma_2 + \nu\,\beta_2 = \displaystyle{\frac{1}{\sqrt{E}} \, \mu_u},\\
\vspace{2mm}
2\mu \, \gamma_1 - \nu \,\beta_1 = \displaystyle{\frac{1}{\sqrt{G}}\, \mu_v},\\
\vspace{2mm}
2\nu\, \gamma_1 - \mu \,\beta_1 = \displaystyle{\frac{1}{\sqrt{G}}\,\nu_v},\\
\vspace{2mm}
2\nu\, \gamma_2 + \mu \,\beta_2 = \displaystyle{\frac{1}{\sqrt{E}}\,\nu_u},\\
\vspace{2mm} \gamma_2\,\beta_2 - \gamma_1\,\beta_1 - 2\nu \mu =
\displaystyle{\frac{1}{\sqrt{E}}\,(\beta_2)_u -
\frac{1}{\sqrt{G}}\,(\beta_1)_v},\\
\vspace{2mm} \gamma_1^2 + \gamma_2^2 - (\nu^2 + \mu^2) =
\displaystyle{\frac{1}{\sqrt{E}}\,(\gamma_2)_u +
\frac{1}{\sqrt{G}}\,(\gamma_1)_v}.
\end{array} \leqno{(5.1)}$$

So, if $\gamma_1 \neq 0$ and $\gamma_2 \neq 0$, then from (5.1) we obtain
$$\sqrt{E} = \displaystyle{\frac{\left(\ln |\mu^2 - \nu^2| \right)_u}{4 \gamma_2}};
\qquad \sqrt{G} = \displaystyle{\frac{\left(\ln |\mu^2 - \nu^2| \right)_v}{4 \gamma_1}}.$$

Following the scheme in \cite{GMih} we give the following definition:

A minimal surface of general type is said to be \emph{strongly regular}
if $\gamma_1 \gamma_2 \neq 0$.
\vskip 2mm
We shall prove the fundamental theorem of Bonnet type for strongly regular minimal
surfaces of general type:

\begin{thm}\label{T:Fundamental Theorem}
Let $ \nu,\, \mu, \, \gamma_1, \, \gamma_2$ be smooth functions,
defined in a domain ${\mathcal D}, \,\, {\mathcal D} \subset
{\R}^2$, and satisfying the conditions
$$\begin{array}{l}
\mu > 0, \qquad \qquad \qquad \qquad \quad \nu \gamma_1 \gamma_2 \neq 0,\\
[3mm]
\gamma_2 \,\left(\ln |\mu^2 - \nu^2| \right)_u >0,
\qquad \,\,
\gamma_1\,\left(\ln |\mu^2 - \nu^2|\right)_v >0,\\
[3mm]
\gamma_1 \sqrt{E} \sqrt{G} = - (\sqrt{E})_v, \qquad
\gamma_2 \sqrt{E} \sqrt{G} = - (\sqrt{G})_u,\\
[3mm]
 \gamma_2\,\beta_2 - \gamma_1\,\beta_1 - 2 \nu \mu =
\displaystyle{\frac{1}{\sqrt{E}}\,(\beta_2)_u -
\frac{1}{\sqrt{G}}\,(\beta_1)_v},\\
[3mm]
\gamma_1^2 + \gamma_2^2 - (\nu^2 + \mu^2) =
\displaystyle{\frac{1}{\sqrt{E}}\,(\gamma_2)_u + \frac{1}{\sqrt{G}}\,(\gamma_1)_v},
\end{array}\leqno{(5.2)}$$
where
\vskip 1mm
\centerline{$\sqrt{E} = \displaystyle{\frac{\left(\ln |\mu^2 - \nu^2|\right)_u}
{4 \gamma_2}}$, $\sqrt{G} = \displaystyle{\frac{\left(\ln |\mu^2 - \nu^2|\right)_v}
{4 \gamma_1}}$, $\beta_1 = 2\gamma_1\, \displaystyle{\frac{\mu \nu_v - \nu \mu_v}
{\nu \nu_v - \mu \mu_v}}$, $\beta_2 = - 2\gamma_2 \,
\displaystyle{\frac{\mu \nu_u - \nu \mu_u}{\nu \nu_u - \mu \mu_u}}$.}
\vskip 1mm
If $\{x_0, \, y_0, \, (n_1)_0,\, (n_2)_0\}$
is a positive oriented orthonormal frame at a point $p_0 \in \R^4$, then there
exists a subdomain ${\mathcal D}_0 \subset {\mathcal D}$ and a unique
strongly regular general minimal surface $M^2: z = z(u,v), \,\,
(u,v) \in {\mathcal D}_0$, passing through $p_0$, such that the invariants of
$M^2$ are the given functions $\nu, \, \mu, \, \gamma_1, \, \gamma_2$, and
$\{x_0, \, y_0, \, (n_1)_0,\, (n_2)_0\}$ is the moving frame of $M^2$ at
the point $p_0$.
\end{thm}
\vskip 2mm
\noindent
\emph{Proof:} We consider the following system of partial differential equations
for the unknown vector functions $x = x(u,v), \, y = y(u,v), \,n_1 = n_1(u,v),
\,n_2 = n_2(u,v)$ in $\R^4$:
$$\begin{array}{ll}
\vspace{2mm}
x_u = \sqrt{E}\, \gamma_1\, y + \sqrt{E}\, \nu \, n_1, & \qquad x_v = - \sqrt{G}\,
\gamma_2\, y + \sqrt{G}\, \mu\, n_2, \\
\vspace{2mm}
y_u = - \sqrt{E}\, \gamma_1\, x + \sqrt{E}\, \mu\, n_,  & \qquad y_v = \sqrt{G}\,
\gamma_2\, x - \sqrt{G}\, \nu\, n_1, \\
\vspace{2mm}
(n_1)_u = - \sqrt{E}\, \nu\, x  + \sqrt{E}\, \beta_1\, n_2,  & \qquad (n_1)_v =
\sqrt{G}\, \nu \, y + \sqrt{G}\, \beta_2\, n_2, \\
\vspace{2mm}
(n_2)_u = - \sqrt{E}\, \mu\, y - \sqrt{E}\, \beta_1\, n_1,  &
\qquad (n_2)_v = - \sqrt{G}\, \mu\, x - \sqrt{G}\, \beta_2\, n_1.
\end{array}\leqno{(5.3)}$$
We denote
$$Z =
\left(%
\begin{array}{c}
  x \\
  y \\
  n_1 \\
  n_2 \\
\end{array}%
\right), \quad
A = \sqrt{E} \left(%
\begin{array}{cccc}
  0 & \gamma_1 & \nu & 0 \\
  -\gamma_1 & 0 & 0 & \mu \\
  -\nu & 0 & 0 & \beta_1 \\
  0 & -\mu & -\beta_1 & 0 \\
\end{array}%
\right), \quad B = \sqrt{G}
\left(%
\begin{array}{cccc}
  0 & -\gamma_2 & 0 & \mu \\
  \gamma_2 & 0 & - \nu & 0 \\
  0 & \nu & 0 & \beta_2 \\
  -\mu & 0 & -\beta_2 & 0 \\
\end{array}%
\right).$$ Then the system (5.3) can be rewritten in the matrix-form:
$$ \begin{array}{l}
\vspace{2mm}
Z_u = A\,Z,\\
\vspace{2mm}
Z_v = B\,Z.
\end{array}\leqno{(5.4)}$$
The integrability conditions of (5.4) are
$$Z_{uv} = Z_{vu},$$
i.e.
$$\displaystyle{\frac{\partial a_i^k}{\partial v} - \frac{\partial b_i^k}{\partial u} +
\sum_{j=1}^{4}(a_i^j\,b_j^k - b_i^j\,a_j^k) = 0, \quad i,k = 1, \dots, 4,} \leqno(5.5)$$
where $a_i^j$ and $b_i^j$ are the elements of the matrices $A$ and $B$. Using (5.2) we
obtain that the equalities (5.5) are fulfilled. Hence, there exists a subset
${\mathcal D}_1 \subset {\mathcal D}$ and unique vector functions $x
= x(u,v), \, y = y(u,v), \,n_1 = n_1(u,v), \,n_2 = n_2(u,v), \,\,
(u,v) \in {\mathcal D}_1$, which satisfy the system (5.3) and the conditions
$$x(u_0,v_0) = x_0, \quad y(u_0,v_0) = y_0, \quad n_1(u_0,v_0)
= (n_1)_0, \quad n_2(u_0,v_0) = (n_2)_0.$$

We shall prove that $\{x(u,v), \, y(u,v), \,n_1(u,v), \,n_2(u,v)\}$
is a positive oriented orthonormal frame field for each $(u,v) \in {\mathcal D}_1$.
Let us consider the following functions:
$$\begin{array}{llll}
\vspace{2mm}
\varphi_1 = x^2 - 1, & \qquad \varphi_2 = y^2 - 1, & \qquad \varphi_3 = n_1^2 - 1,
& \qquad \varphi_4 = n_2^2 - 1,\\
\vspace{2mm}
\varphi_5 = x\,y,  & \qquad \varphi_6 = x\,n_1, & \qquad \varphi_7 = x\,n_2,
& \qquad \varphi_8 = y\,n_1,\\
\vspace{2mm}
& \qquad \varphi_9 = y\,n_2, & \qquad \varphi_{10} = n_1\,n_2, &\\
\end{array}$$
defined for each $(u,v) \in {\mathcal D}_1$. Using that
$x(u,v), \, y(u,v), \,n_1(u,v), \,n_2(u,v)$ satisfy (5.3), we obtain  the system
$$\begin{array}{lll}
\vspace{2mm}
\displaystyle{\frac{\partial \varphi_i}{\partial u} = \alpha_i^j \, \varphi_j},\\
\vspace{2mm}
\displaystyle{\frac{\partial \varphi_i}{\partial v} = \beta_i^j \, \varphi_j},
\end{array}
\qquad i = 1, \dots, 10, \leqno{(5.6)}$$
where $\alpha_i^j, \beta_i^j, \,\, i,j = 1, \dots, 10$ are functions of
$(u,v) \in {\mathcal D}_1$. The system (5.6) is a linear system of partial differential
equations for the functions
$\varphi_i(u,v), \,\,i = 1, \dots, 10, \,\,(u,v) \in {\mathcal D}_1$, satisfying
$\varphi_i(u_0,v_0) = 0, \,\,i = 1, \dots, 10$. Hence
$\varphi_i(u,v) = 0, \,\,i = 1, \dots, 10$ for each $(u,v) \in {\mathcal D}_1$.
Consequently, the vector functions $x(u,v), \, y(u,v), \,n_1(u,v), \,n_2(u,v)$
form an orthonormal frame field for each $(u,v) \in {\mathcal D}_1$.

Now, let us consider the  system
$$\begin{array}{l}
\vspace{2mm}
z_u = \sqrt{E}\, x,\\
\vspace{2mm} z_v = \sqrt{G}\, y
\end{array}\leqno{(5.7)}$$
of partial differential equations for the vector function $z(u,v)$. Using (5.2)
and (5.3) we get that the integrability conditions $z_{uv} = z_{vu}$ of (5.7) are
fulfilled. Hence,  there exist a subset ${\mathcal D}_0 \subset {\mathcal D}_1$ and
a unique vector function $z = z(u,v)$, defined for $(u,v) \in {\mathcal D}_0$ and
satisfying $z(u_0, v_0) = p_0$.

Consequently, the surface $M^2: z = z(u,v), \,\, (u,v) \in {\mathcal D}_0$ satisfies
the assertion of the theorem. \qed

\section{Canonical parameters on strongly regular minimal surfaces of general type}

Let $M^2$ be a strongly regular minimal surface of general type parameterized by
semi-canonical parameters. In this section we shall define canonical parameters on $M^2$.

\begin{lem}\label{L:Canonocal}
Let $M^2$ be a strongly regular minimal surface of general type in $\R^4$,
parameterized by semi-canonical parameters. Then the function
$$E \sqrt{|\mu^2 - \nu^2|}$$
does not depend on $v$, while the function
$$G \sqrt{|\mu^2 - \nu^2|}$$
does not depend on $u$.
\end{lem}
\vskip 2mm
\noindent
\emph{Proof:} From equalities (5.1) we get
$$\gamma_1 = \ds{\frac{1}{4 \sqrt{G}}\left(\ln |\mu^2 - \nu^2|\right)_v}; \qquad
\gamma_2 = \ds{\frac{1}{4 \sqrt{E}}\left(\ln |\mu^2 - \nu^2|\right)_u}. \leqno{(6.1)}$$
On the other hand,
$\gamma_1 = \ds{- \frac{1}{\sqrt{G}} \left(\ln\sqrt{E}\right)_v}, \,\, \gamma_2
=  \ds{- \frac{1}{\sqrt{E}} \left(\ln\sqrt{G}\right)_u}$. Hence,
$$\ds{\frac{\partial}{\partial v}\left(\ln \left(E^2 |\mu^2 - \nu^2|\right)\right)} = 0;
\qquad \ds{\frac{\partial}{\partial u}\left(\ln \left(G^2 |\mu^2 -
\nu^2|\right)\right)} = 0,$$ which imply that $E \sqrt{|\mu^2 - \nu^2|}$ does not
depend on $v$, and $G \sqrt{|\mu^2 - \nu^2|}$ does not depend on $u$. \qed
\vskip 2mm
The equalities (6.1) imply the following equivalences:
$$\gamma_1 \neq 0 \qquad \Longleftrightarrow \qquad |\mu^2 - \nu^2|_v \neq 0;$$
$$\gamma_2 \neq 0 \qquad \Longleftrightarrow \qquad |\mu^2 - \nu^2|_u \neq 0.$$
So, the strongly regular minimal surfaces of general type are characterized
by the condition
$$|\mu^2 - \nu^2|_v |\mu^2 - \nu^2|_u \neq 0.$$

Next we introduce canonical parameters on a strongly regular minimal surface of
general type.

\begin{defn}\label{D:Canonocal parameters}
Let $M^2$ be a strongly regular minimal surface of general type parameterized by
semi-canonical parameters $(u,v)$. The parameters $(u,v)$ are said to be
\textit{canonical}, if
$$E \sqrt{|\mu^2 - \nu^2|} = 1; \qquad G \sqrt{|\mu^2 - \nu^2|} = 1.$$
\end{defn}

\begin{thm}\label{T:Canonocal parameters}
Any strongly regular minimal surface of general type locally admits canonical parameters.
\end{thm}
\vskip 2mm
\noindent
\emph{Proof:} Let $M^2$ be a strongly regular minimal surface of general type
parameterized by semi-canonical parameters $(u,v)$. According to Lemma \ref{L:Canonocal}
it follows that there exist functions $\varphi = \varphi(u) >0 $ and
$\psi = \psi(v) >0$, such that
$$E \sqrt{|\mu^2 - \nu^2|} = \varphi(u); \qquad G \sqrt{|\mu^2 - \nu^2|} = \psi(v).$$
Under the following change of the parameters:
$$\begin{array}{l}
\vspace{2mm}
\overline{u} = \ds{\int_{u_0}^u \sqrt{\varphi(u)}\, du}
+ \overline{u}_0, \quad \overline{u}_0 = const\\
\vspace{2mm}
\overline{v} = \ds{\int_{v_0}^v  \sqrt{\psi(v)}\, dv +
\overline{v}_0}, \quad \overline{v}_0 = const
\end{array}$$
we obtain
$$\overline{E} = \ds{\frac{1}{ \sqrt{|\mu^2 - \nu^2|}}}; \qquad \overline{F} = 0;
\qquad \overline{G} = \ds{\frac{1}{ \sqrt{|\mu^2 - \nu^2|}}},$$ i.e. the parameters
$(\overline{u},\overline{v})$ are canonical. \qed
\vskip 2mm
From now on we consider a strongly regular minimal surface of general type
$M^2: z = z(u,v), \,\, (u,v) \in {\mathcal D}$, parameterized by canonical
parameters, i.e.
$$E = \ds{\frac{1}{ \sqrt{|\mu^2 - \nu^2|}}}; \qquad
F = 0; \qquad G = \ds{\frac{1}{ \sqrt{|\mu^2 - \nu^2|}}}.$$
Then the functions $\gamma_1$ and $\gamma_2$ in the derivative formulas of
$M^2$ are expressed as follows:
$$\gamma_1 = \left(|\mu^2 - \nu^2|^{\frac{1}{4}}\right)_v; \qquad \gamma_2
= \left(|\mu^2 - \nu^2|^{\frac{1}{4}}\right)_u.$$
Using equalities (5.1) we find that
$$\beta_1 = - \ds{|\mu^2 - \nu^2|^{\frac{1}{4}}
\left(\ln\sqrt{\left|\frac{\mu + \nu}{\mu - \nu}\right|}\right)_v};
\qquad \beta_2 = \ds{|\mu^2 - \nu^2|^{\frac{1}{4}}
\left(\ln\sqrt{\left|\frac{\mu + \nu}{\mu - \nu}\right|}\right)_u}.$$
Then the last two equalities of (5.1) take the form:
$$\begin{array}{l}
\vspace{2mm}
\ds{\frac{1}{4}}\sqrt{|\mu^2 - \nu^2|}\, \Delta \ln |\mu^2 - \nu^2|
+ \nu^2 + \mu^2 = 0;\\
\vspace{2mm}
\ds{\frac{1}{2}}\sqrt{|\mu^2 - \nu^2|} \, \Delta \ln
\left|\frac{\mu + \nu}{\mu - \nu}\right| + 2\nu \mu = 0,
\end{array} \leqno{(6.2)}$$
where $\Delta$ is the Laplace operator.

The invariants $\mu$ and $\nu$ can be expressed by the Gauss
curvature $K$ and the normal curvature $\varkappa$ as follows:
$$|\mu^2 - \nu^2| = \sqrt{K^2 - \varkappa^2}; \qquad \ds{\left|\frac{\mu + \nu}
{\mu - \nu}\right| = \sqrt{\frac{K - \varkappa}{K + \varkappa}}}.$$

So the equalities (6.2) can be rewritten as
$$\begin{array}{l}
\vspace{2mm}
\ds{\frac{1}{8}}(K^2 - \varkappa^2)^{\frac{1}{4}}\, \Delta \ln (K^2 - \varkappa^2) - K = 0;\\
\vspace{2mm}
\ds{\frac{1}{4}}(K^2 - \varkappa^2)^{\frac{1}{4}} \,
\Delta \ln \frac{K - \varkappa}{K + \varkappa} + \varkappa = 0.
\end{array} \leqno{(6.3)}$$

The fundamental Theorem \ref{T:Fundamental Theorem} in canonical parameters states as follows:

\begin{thm}\label{T:Fundamental Theorem - Canonical}
Let $\nu(u,v)$ and $\mu(u,v)$ be two smooth functions, defined in a domain
${\mathcal D}, \,\, {\mathcal D} \subset {\R}^2$, and satisfying the conditions
$$\begin{array}{l}
\vspace{2mm}
\mu > 0, \qquad \nu \neq 0, \qquad |\mu^2 - \nu^2|_v |\mu^2 - \nu^2|_u \neq 0;\\
\vspace{2mm}
\ds{\frac{1}{4}}\sqrt{|\mu^2 - \nu^2|}\, \Delta \ln |\mu^2 - \nu^2| + \nu^2 + \mu^2 = 0;\\
\vspace{2mm} \ds{\frac{1}{2}}\sqrt{|\mu^2 - \nu^2|} \, \Delta \ln
\left|\frac{\mu + \nu}{\mu - \nu}\right| + 2\nu \mu = 0.
\end{array} $$
Then there exists a unique (up to a motion) strongly regular minimal surface of
general type
$M^2: z = z(u,v), \,\, (u,v) \in {\mathcal D}_0,
\,\, {\mathcal D}_0 \subset {\mathcal D}$, with geometric invariants
\, $\nu(u,v)$ and $\mu(u,v)$. Furthermore $(u,v)$ are canonical parameters of $M^2$.
\end{thm}

\section{Minimal surfaces of general type which are not strongly regular}

In this section we consider non strongly regular minimal surfaces of general type,
satisfying the condition
$$\gamma_1 = 0. \leqno{(7.1)}$$

Equalities (5.1) under the assumption $\gamma_1 = 0$ imply that
$$\begin{array}{l}
\vspace{2mm}
\nu \,\beta_1 = - \displaystyle{\frac{1}{\sqrt{G}}\, \mu_v},\\
\vspace{2mm}
\mu \,\beta_1 = - \displaystyle{\frac{1}{\sqrt{G}}\,\nu_v},\\
\vspace{2mm}
2 \mu \, \gamma_2 + \nu\,\beta_2 = \displaystyle{\frac{1}{\sqrt{E}} \, \mu_u},\\
\vspace{2mm}
2\nu\, \gamma_2 + \mu \,\beta_2 = \displaystyle{\frac{1}{\sqrt{E}}\,\nu_u},\\
\vspace{2mm} \gamma_2\,\beta_2 - 2\nu \mu =
\displaystyle{\frac{1}{\sqrt{E}}\,(\beta_2)_u -
\frac{1}{\sqrt{G}}\,(\beta_1)_v},\\
\vspace{2mm} \gamma_2^2 - (\nu^2 + \mu^2) =
\displaystyle{\frac{1}{\sqrt{E}}\,(\gamma_2)_u}.
\end{array} \leqno{(7.2)}$$

Having in mind that $\gamma_1= -y(\ln\sqrt{E})$, and (7.1), we get
$E_v = 0$, i.e. $E = E(u)$. From the first two equalities of (7.2)
we obtain $(\mu^2 - \nu^2)_v = 0$. The third and the fourth equality
of (7.2) imply
$$\gamma_2 = \ds{\frac{1}{4 \sqrt{E}}\left(\ln |\mu^2 - \nu^2|\right)_u}, \qquad
\beta_2 = \ds{{\frac{1}{\sqrt{E}}} \left(\ln\sqrt{\left|\frac{\mu +
\nu}{\mu - \nu}\right|}\right)_u}. \leqno{(7.3)}$$
Hence, $\gamma_2 = \gamma_2(u)$. Then from the last equality of (7.2) we get that
$(\mu^2 + \nu^2)_v = 0$. Consequently, $\mu = \mu(u), \,\, \nu = \nu(u)$. So, we
obtain that $\beta_1 = 0$ and $\beta_2 = \beta_2(u)$.

With respect to the geometric frame field $\{x, y, n_1, n_2\}$ the
derivative formulas of $M^2$ look like:
$$\begin{array}{ll}
\vspace{2mm} \nabla'_xx = \quad \quad \quad \qquad  \,\nu\, n_1, &
\qquad
\nabla'_x n_1 = - \nu\,x,\\
\vspace{2mm}
\nabla'_xy = \qquad \quad \quad \quad \; \qquad \quad + \mu\, n_2, & \qquad
\nabla'_y n_1 = \quad \qquad \nu \,y \quad \quad \quad + \beta_2\, n_2,\\
\vspace{2mm}
\nabla'_yx = \quad \quad \; - \gamma_2\,y \; \qquad \quad + \mu \, n_2,  & \qquad
\nabla'_x n_2= \qquad - \,\, \mu \,y,\\
\vspace{2mm}
\nabla'_yy = \;\;\gamma_2\,x \quad \quad \,\,\, - \,\nu \, n_1, & \qquad
\nabla'_y n_2 = - \mu \,x \qquad \,\, - \beta_2\,n_1,
\end{array} \leqno{(7.4)}$$
where $ \mu > 0,\,\, \nu \neq 0$, $\mu^2 \neq \nu ^2$, $\gamma_2 \neq0$, $\beta_2 \neq 0$.

\vskip 2mm
We shall define canonical parameters on $M^2$.

\begin{lem}\label{L:Canonocal-1}
Let $M^2$ be a minimal surface of general type parameterized by
semi-canonical parameters, and satisfying the condition $\gamma_1 = 0$. Then
the function $G \sqrt{|\mu^2 - \nu^2|}$ does not depend on $u$.
\end{lem}
\noindent
\emph{Proof:} Using that $\gamma_2= - \ds{\frac{1}{\sqrt{E}}(\ln\sqrt{G})_u}$
and (7.3) we obtain
$\ds{\frac{\partial}{\partial u}\left(\ln \left(G^2 |\mu^2 - \nu^2|\right)\right)} = 0,$
which implies that $G \sqrt{|\mu^2 - \nu^2|}$ does not depend on $u$. \qed
\vskip 2mm
\begin{defn}\label{D:Canonocal parameters-1}
Let $M^2$ be a minimal surface of general type parameterized by semi-canonical
parameters $(u,v)$, and satisfying the condition $\gamma_1 = 0$.
The parameters $(u,v)$ are said to be \textit{canonical}, if
$$E = 1, \qquad G = \ds{\frac{1}{\sqrt{|\mu^2 - \nu^2|}}}.$$
\end{defn}

\begin{thm}\label{T:Canonocal parameters-1}
Any minimal surface of general type satisfying the condition $\gamma_1 =0$ locally
admits canonical parameters.
\end{thm}
\vskip 2mm
\noindent
\emph{Proof:} Let $M^2$ be a non-strongly regular minimal surface of general type,
parameterized by semi-canonical parameters $(u,v)$, and satisfying $\gamma_1 = 0$.
According to Lemma \ref{L:Canonocal-1} it follows that there exist functions
$\varphi = \varphi(u) >0 $ and $\psi = \psi(v) >0$, such that
$$E = \varphi(u), \qquad G \sqrt{|\mu^2 - \nu^2|} = \psi(v).$$
Under the following change of the parameters:
$$\begin{array}{l}
\vspace{2mm}
\overline{u} = \ds{\int_{u_0}^u \sqrt{\varphi(u)}\, du} + \overline{u}_0,
\quad \overline{u}_0 = const\\
\vspace{2mm} \overline{v} = \ds{\int_{v_0}^v  \sqrt{\psi(v)}\, dv +
\overline{v}_0}, \quad \overline{v}_0 = const
\end{array}$$
we get
$$\overline{E} = 1, \qquad \overline{F} = 0, \qquad \overline{G} =
\ds{\frac{1}{ \sqrt{|\mu^2 -
\nu^2|}}},$$
i.e. the parameters $(\overline{u},\overline{v})$ are canonical. \qed
\vskip 2mm
From now on we consider a minimal surface $M^2$ of general type satisfying the
condition $\gamma_1=0$, parameterized by canonical parameters, i.e.
$$E =1, \qquad F = 0, \qquad G = \ds{\frac{1}{ \sqrt{|\mu^2 - \nu^2|}}}.$$ Then
the functions $\gamma_1$, $\gamma_2$, $\beta_1$,  $\beta_2$ in the
derivative formulas (7.4) of $M^2$ are expressed as follows:
$$\begin{array}{ll}
\vspace{2mm}
\gamma_1 = 0, \qquad & \gamma_2 = \ds{\frac{1}{4}\left(\ln |\mu^2 - \nu^2|\right)_u};\\
\vspace{2mm}
\beta_1 = 0, \qquad & \beta_2 =
\ds{\frac{1}{2}\left(\ln \left|\frac{\mu + \nu}{\mu - \nu}\right|
\right)_u}.
\end{array}$$
Hence, the last two equalities of (7.2) take the form:
$$\begin{array}{l}
\vspace{2mm}
\ds{\frac{1}{4}\left(\ln |\mu^2 - \nu^2|\right)_{uu}
- \frac{1}{16}\left(\ln |\mu^2 - \nu^2|\right)_u^2  + \nu^2 + \mu^2 = 0},\\
\vspace{2mm}
\ds{\frac{1}{2}\left(\ln \left|\frac{\mu + \nu}{\mu - \nu}\right|\right)_{uu}
- \frac{1}{8}\left(\ln |\mu^2 - \nu^2|\right)_u \left(\ln \left|\frac{\mu + \nu}
{\mu - \nu}\right|\right)_u + 2\nu \mu = 0}.
\end{array}  \leqno{(7.5)}$$

\begin{thm}\label{T:Main Theorem - Canonical-1}
Let $\nu(u)$ and $\mu(u)$ be two smooth functions, defined in $J
\subset \R$, and satisfying the conditions
$$\begin{array}{l}
\vspace{2mm} \mu > 0; \qquad \nu \neq 0; \qquad |\mu^2 - \nu^2|_u
\neq 0; \qquad
\ds{\left|\frac{\mu + \nu}{\mu - \nu}\right|_u \neq 0};\\
\vspace{2mm}
\ds{\frac{1}{4}\left(\ln |\mu^2 - \nu^2|\right)_{uu}
- \frac{1}{16}\left(\ln |\mu^2 - \nu^2|\right)_u^2  + \nu^2 + \mu^2 = 0};\\
\vspace{2mm}
\ds{\frac{1}{2}\left(\ln \left|\frac{\mu + \nu}{\mu - \nu}\right|\right)_{uu}
- \frac{1}{8}\left(\ln |\mu^2 - \nu^2|\right)_u \left(\ln \left|\frac{\mu + \nu}
{\mu - \nu} \right|\right)_u + 2\nu \mu = 0}.
\end{array} $$
Then there exists a unique (up to a motion) minimal surface of general type
$M^2: z = z(u,v), \,\, (u,v) \in {\mathcal D}_0, \,\,
{\mathcal D}_0 =J\times J' \; (J'\subset \R)$, satisfying the condition
$\gamma_1=0$  with geometric invariants $\nu(u)$ and $\mu(u)$. Furthermore $(u,v)$
are canonical parameters of $M^2$.
\end{thm}
\vskip 2mm
\noindent
\emph{Proof:}
Let us denote
$$E =1; \quad G = \ds{\frac{1}{ \sqrt{|\mu^2 - \nu^2|}}}; \quad
\gamma_2 = \ds{\frac{1}{4}\left(\ln |\mu^2 - \nu^2|\right)_u}; \quad
\beta_2 = \ds{\frac{1}{2}\left(\ln \left|\frac{\mu + \nu}{\mu - \nu}\right| \right)_u}.$$
We consider the following system of partial
differential equations for the unknown vector functions $x = x(u,v),
\, y = y(u,v), \,n_1 = n_1(u,v), \,n_2 = n_2(u,v)$ in $\R^4$:
$$\begin{array}{ll}
\vspace{2mm}
x_u = \nu\, n_1, & \qquad x_v = - \sqrt{G}\,
\gamma_2\, y + \sqrt{G}\, \mu\, n_2, \\
\vspace{2mm}
y_u = \mu\, n_2,  & \qquad y_v = \sqrt{G}\,
\gamma_2\, x - \sqrt{G}\, \nu\, n_1, \\
\vspace{2mm}
(n_1)_u = - \nu\, x,   & \qquad (n_1)_v =
\sqrt{G}\, \nu \, y + \sqrt{G}\, \beta_2\, n_2, \\
\vspace{2mm}
(n_2)_u = -  \mu\, y,  &
\qquad (n_2)_v = - \sqrt{G}\, \mu\, x - \sqrt{G}\, \beta_2\, n_1.
\end{array}\leqno{(7.6)}$$
The integrability conditions of system (7.6) are equalities (7.5).
Further we follow the scheme of the proof of Theorem \ref{T:Fundamental Theorem}.
\qed
\vskip 2mm
We shall give a geometric description of the non strongly regular
minimal surfaces of general type, satisfying the condition $\gamma_1=0$.

\begin{prop}\label{P:On a class}
Let $M^2$ be a minimal surface of general type, satisfying the condition
$\gamma_1 = 0$. Then each parametric curve $c_v: z(v) = z(u_0,v)$ is
a curve with constant curvatures, and each parametric curve $c_u: z(u) = z(u,v_0)$
is a plane curve, lying in the $2$-dimensional space spanned by the second and
the fourth vector field in the Frenet frame field of $c_v$.
\end{prop}
\vskip 2mm
\noindent
\emph{Proof:} Let $c_u: z(u) = z(u,v_0)$ be a parametric
curve of $M^2$. Then $z'(u) = z_u(u,v_0)$ and the tangent
vector $t_u$ of $c_u$ is $t_u = \ds{\frac{z_u}{\sqrt{E}}} = x(u,v_0)$. Using (7.4)
we obtain $t_u' = \nu\, n_1$, which implies that the curvature of $c_u$ is
$\varkappa_u(u) = |\nu(u)|$, and the principal normal $n_u$ of $c_u$ is collinear
with $n_1$. Then from (7.4) we get $n_u' = -\varkappa_u t_u$, which implies that the
parametric curve $c_u$ is a plane curve, lying in the 2-dimensional
space $\span\{x,n_1\}$.

\vskip 2mm Now, let $c_v: z(v) = z(u_0,v)$ be a parametric curve of
$M^2$. Then $z'(v) = z_v(u_0,v)$ and the tangent vector $t_v$
of $c_v$ is $t_v = \ds{\frac{z_v}{\sqrt{G}}} = y(u_0,v)$. Using
(7.4) we find the Frenet frame field of $c_v$  and we calculate the
curvatures of $c_v$. The Frenet frame field $\{t_v, n_v, b_v, b_{1\,v}\}$ of $c_v$
is determined by
$$t_v = y, \quad  n_v = \ds{\frac{\gamma_2}{\sqrt{\nu^2
+ \gamma_2^2}}\,x - \frac{\nu}{\sqrt{\nu^2 + \gamma_2^2}}\,n_1},
\quad  b_v = n_2, \quad  b_{1\,v} = \ds{\frac{\nu}{\sqrt{\nu^2 +
\gamma_2^2}}\,x + \frac{\gamma_2}{\sqrt{\nu^2 + \gamma_2^2}}\,n_1}.$$
Obviously, $\span\{n_v, b_{1\,v}\} = \span\{x,n_1\}$.

The curvatures of $c_v$ are given by
$$\varkappa_v = \varkappa_v(u_0) = \sqrt{\nu^2 + \gamma_2^2}, \quad
\tau_v = \tau_v(u_0) = \ds{\frac{|\mu \gamma_2 - \nu \beta_2|}
{\sqrt{\nu^2 + \gamma_2^2}}}, \quad \sigma_v = \sigma_v(u_0) = -
\ds{\frac{\mu \nu + \gamma_2 \beta_2} {\sqrt{\nu^2 + \gamma_2^2}}}.$$
Since the curvatures of $c_v$ do not depend on the parameter $v$, the parametric
curve $c_v$ is a curve with constant curvatures. \qed
\vskip 3mm
The following question arises naturally:
\textit{Does each smooth curve with constant curvatures in $\R^4$ generate a
non strongly regular minimal surface in this geometric way?}

We shall construct a class of surfaces in $\R^4$, generated by a curve with constant
curvatures.
\vskip 2mm
Each curve with constant curvatures in $\R^4$ can be parameterized as follows:
$$c: z(v) = \left( a \cos \alpha v, a \sin \alpha v, b \cos \beta v,
b \sin \beta v \right), \quad  v \in [0; 2\pi),$$
where $a,\,b,\,\alpha, \,\beta$ are constants ($\alpha > 0,\,\,\beta > 0$).
This curve is a generalization of the circular helix in $\E^3$.
Without loss of generality we can assume that $a^2 \alpha^2 + b^2 \beta^2 = 1$,
i.e. $c$ is parameterized by the arc-length. Then the curvatures of $c$ are
expressed by
$$\varkappa = \ds{\sqrt{a^2 \alpha^4 + b^2 \beta^4}}; \quad
\tau = \ds{\frac{ab \alpha \beta (\alpha^2 - \beta^2)}{\sqrt{a^2
\alpha^4 + b^2 \beta^4}}}; \quad \sigma = \ds{\frac{\alpha
\beta}{\sqrt{a^2 \alpha^4 + b^2 \beta^4}}}.$$
In case of $\alpha \neq \beta$ $c$ is a curve in $\R^4$, and in case of
$\alpha = \beta$ $c$ is a circle. We shall consider the general case
$\alpha \neq \beta$.

The Frenet frame field $\{t, n, b, b_1\}$ of $c$ is
$$\begin{array}{l}
\vspace{2mm}
t = \left(- a \alpha \sin \alpha v, a \alpha \cos \alpha v, - b \beta \sin \beta v,
b \beta \cos \beta v \right),\\
\vspace{2mm}
n =\ds{\frac{1}{\varkappa}} \left(- a \alpha^2 \cos \alpha v, - a \alpha^2 \sin \alpha v,
- b \beta^2 \cos \beta v, - b \beta^2 \sin \beta v \right),\\
\vspace{2mm}
b = \left(b \alpha \beta \cos \alpha v, b \alpha \beta \sin \alpha v,
- a \alpha \beta \cos \beta v, - a \alpha \beta \sin \beta v \right),\\
\vspace{2mm}
b_1 =\ds{\frac{1}{\varkappa}} \left(b \beta^2 \cos \alpha v, b \beta^2 \sin \alpha v,
- a \alpha^2 \cos \beta v, - a \alpha^2 \sin \beta v \right).
\end{array}$$

Let $A = A(u)$ and $B = B(u)$ be smooth functions, defined in $J \subset \R$, and satisfying
$$A'\,^2(u)+ B'\,^2(u) > 0, \qquad (\varkappa A - 1)^2 + (\tau A - \sigma B)^2 > 0;
\quad u \in J.$$
We consider the surface $M^2$ generated by the curve $c$ in the following way:
$$M^2: z(u,v) = z(v) + A(u) n(v) + B(u) b_1(v), \quad u \in J, \, v \in [0; 2\pi). \leqno (7.7)$$
The surface $M^2$ is a one-parameter family of plane curves lying in the normal plane
$\span\{n, b_1\}$ of $c$.

Let us denote
$$\begin{array}{ll}
\vspace{2mm}
f(u) = a + \ds{\frac{1}{\varkappa} \left( - a \alpha^2 \, A(u) + b \beta^2\,B(u)\right)},\\
\vspace{2mm}
g(u) = b + \ds{\frac{1}{\varkappa} \left( - b \beta^2 \, A(u) - a \alpha^2 \,B(u)\right)}.
\end{array} \leqno{(7.8)}$$
Then the surface $M^2$ is parameterized by
$$M^2: z(u,v) = \left( f(u) \cos \alpha v, \, f(u) \sin \alpha v, \, g(u) \cos \beta v,
\, g(u) \sin \beta v \right). \leqno{(7.9)}$$

We shall consider a surface $M^2$ given by (7.9), for arbitrary smooth functions $f(u)$ and
$g(u)$, satisfying $\alpha^2 f^2+ \beta^2 g^2 > 0 , \,\, f'\,^2+ g'\,^2 > 0$.

Considering general rotations in $\R^4$, Moore introduces in \cite{M1} general rotational
surfaces. The surface $M^2$, given by (7.9) is a general rotational surface whose meridians
lie in two-dimensional planes.

The tangent space of $M^2$ is spanned by the vector fields
$$\begin{array}{l}
\vspace{2mm}
z_u = \left(f' \cos \alpha v, f' \sin \alpha v, g' \cos \beta v, g' \sin \beta v \right),\\
\vspace{2mm}
z_v = \left( - \alpha f \sin \alpha v, \alpha f \cos \alpha v, - \beta g \sin \beta v,
\beta g \cos \beta v \right).
\end{array}$$
Hence, the coefficients of the first fundamental form are
$$E = f'\,^2(u)+ g'\,^2(u), \quad F = 0, \quad G = \alpha^2 f^2(u)+ \beta^2 g^2(u)$$
and $W = \sqrt{(f'\,^2 + g'\,^2)(\alpha^2 f^2 + \beta^2 g^2)}$. We
consider the following orthonormal tangent frame field
$$\begin{array}{l}
\vspace{2mm}
x = \ds{\frac{1}{\sqrt{f'\,^2 + g'\,^2}}\left(f' \cos \alpha v,
f' \sin \alpha v, g' \cos \beta v, g' \sin \beta v \right)},\\
\vspace{2mm}
y = \ds{\frac{1}{\sqrt{\alpha^2 f^2 + \beta^2 g^2}}\left( - \alpha f \sin \alpha v,
\alpha f \cos \alpha v, - \beta g \sin \beta v, \beta g \cos \beta v \right)}.
\end{array}$$
The second partial derivatives of $z(u,v)$ are expressed as follows
$$\begin{array}{l}
\vspace{2mm}
z_{uu} = \left(f'' \cos \alpha v, f'' \sin \alpha v, g'' \cos \beta v,
g'' \sin \beta v \right),\\
\vspace{2mm} z_{uv} = \left(- \alpha f' \sin \alpha v, \alpha f' \cos \alpha v,
- \beta g' \sin \beta v, \beta g' \cos \beta v \right),\\
\vspace{2mm}
z_{vv} = \left(- \alpha ^2 f \cos \alpha v, - \alpha ^2 f \sin \alpha v,
- \beta ^2 g \cos \beta v, - \beta ^2 g \sin \beta v \right).
\end{array}$$
Now let us consider the following orthonormal normal frame field of $M^2$:
$$\begin{array}{l}
\vspace{2mm}
n_1 = \ds{\frac{1}{\sqrt{f'\,^2 + g'\,^2}}\left(g' \cos \alpha v,
g' \sin \alpha v, - f' \cos \beta v, - f' \sin \beta v \right)},\\
\vspace{2mm}
n_2 = \ds{\frac{1}{\sqrt{\alpha^2 f^2 + \beta^2 g^2}} \left( - \beta g \sin \alpha v,
\beta g \cos \alpha v, \alpha f \sin \beta v, - \alpha f \cos \beta v \right)}.
\end{array}$$
It is easy to verify that $\{x, y, n_1, n_2\}$ is a positive
oriented orthonormal frame field in $\R^4$.

We calculate the functions $c_{ij}^k, \,\, i,j,k = 1,2$:
$$\begin{array}{ll}
\vspace{2mm}
c_{11}^1 = \langle z_{uu}, n_1 \rangle = \ds{\frac{g' f'' - f' g''}
{\sqrt{f'\,^2 + g'\,^2}}}, \quad \quad & c_{11}^2 = \langle z_{uu}, n_2 \rangle = 0,\\
\vspace{2mm}
c_{12}^1 = \langle z_{uv}, n_1 \rangle = 0, \quad \quad
& c_{12}^2 = \langle z_{uv}, n_2 \rangle =
\ds{\frac{\alpha \beta (g f' - f g')}{\sqrt{\alpha^2 f^2 + \beta^2 g^2}}},\\
\vspace{2mm}
c_{22}^1 = \langle z_{vv}, n_1 \rangle =
\ds{\frac{\beta^2 g f' - \alpha^2 f g'}{\sqrt{f'\,^2 + g'\,^2}}},
\quad \quad & c_{22}^2 = \langle z_{vv}, n_2 \rangle = 0,
\end{array}$$
where $\langle \,,\, \rangle$ is the standard scalar product in $\R^4$. Therefore
the coefficients $L$, $M$ and $N$ of the second
fundamental form of $M^2$ are expressed as follows:
$$L = \ds{\frac{2 \alpha \beta (g f' - f g') (g' f'' - f' g'')}
{(\alpha^2 f^2 + \beta^2 g^2) (f'\,^2 + g'\,^2)}}, \qquad M = 0,
\qquad N = \ds{\frac{- 2\alpha \beta (g f' - f g') (\beta^2 g f' -
\alpha^2 f g')}{(\alpha^2 f^2 + \beta^2 g^2) (f'\,^2 + g'\,^2)}}.$$

Consequently, the invariants $k$ and $\varkappa$ of $M^2$ are given by:
$$k = \ds{\frac{- 4 \alpha^2 \beta^2 (g f' - f g')^2 (g' f'' - f' g'')
(\beta^2 g f' - \alpha^2 f g')}{(\alpha^2 f^2 + \beta^2 g^2)^3
(f'\,^2 + g'\,^2)^3}};$$
$$\varkappa =  \ds{\frac{\alpha \beta (g f' - f g')}{(\alpha^2 f^2 +
\beta^2 g^2)^2 (f'\,^2 + g'\,^2)^2} \, [(\alpha^2 f^2 + \beta^2
g^2)(g' f'' - f' g'') - (f'\,^2 + g'\,^2) (\beta^2 g f' - \alpha^2 f g') ]}.$$

Hence, $M^2$ is a minimal surface ($\varkappa^2 - k = 0$) if and
only if
$$\ds{\frac{g' f'' - f' g''}{f'\,^2 + g'\,^2}} = \ds{\frac{\alpha^2 f g'
 - \beta^2 g f'}{\alpha^2 f^2 + \beta^2 g^2}}. \leqno(7.10)$$

Now let $f(u)$ and $g(u)$ be the functions, defined by (7.8). Then
we calculate
$$\begin{array}{l}
\vspace{2mm}
g' f'' - f' g'' = A'' B' - A' B'', \qquad \qquad f'\,^2 + g'\,^2 = A'\,^2 + B'\,^2, \\
\vspace{2mm} \alpha^2 f g' - \beta^2 g f' = (\varkappa A -
1)\varkappa B' + (\tau A - \sigma B) (\sigma A' + \tau B'),\\
\vspace{2mm} \alpha^2 f^2 + \beta^2 g^2 = (\varkappa A - 1)^2 + (\tau A - \sigma B)^2.
\end{array}$$

Consequently, the equality (7.10) takes the form
$$\ds{\frac{A'' B' - A' B''}{A'\,^2 + B'\,^2}} =
\ds{\frac{(\varkappa A - 1)\varkappa B' + (\tau A - \sigma B)
(\sigma A' + \tau B')} {(\varkappa A - 1)^2 + (\tau A - \sigma B)^2}}. \leqno(7.11)$$

Thus we obtain the following result:
\begin{prop}\label {P:Examples}
The surface $M^2$, generated by the curve $c$ with constant curvatures $\varkappa$,
$\tau$ and $\sigma$ by the formula $(7.7)$, is minimal if and only if the equality
$(7.11)$ holds.
\end{prop}

Any minimal surface given by (7.7) satisfies the condition $\gamma_1=0$.

\end{document}